# Investigation of the Estimation Accuracy of 5 Different Numerical ODE Solvers on 3 Case Studies

Hamidreza Moradi[1,*] , Hamideh Hossei[2]

[1]Department of Mechanical Engineering and Engineering Science, The University of North Carolina at Charlotte, Charlotte, North Carolina, USA

[2]Department of Infrastructure and Environmental System, The University of North Carolina at Charlotte, North Carolina, USA

* Corresponding Author: hmoradi@uncc.edu

## Abstract

Numerical ordinary differential equation (ODE) solvers are indispensable tools in various engineering domains, enabling the simulation and analysis of dynamic systems. In this work, we utilize 5 different numerical ODE solvers namely: Euler's method, Heun's method, Midpoint Method, Runge-kutta 4$^{th}$ order and ODE45 method to discover the answer of three well-known case studies and compare their results by calculation of relative errors. To check for the validity of the estimations, the experimental data of previous literature have been compared with the data in this paper which shows a good accordance. We observe that for each of the case studies based on the behavior of the model, the estimation accuracy of the solvers is different. For the logistic population change as the first case study, the results of all solvers are so close to each other that only their solution cost can be considered for their superiority. For temperature change of a building as the second case study we see that in some especial areas the accuracy of the solvers is different and in general Midpoint ODE solver shows better results. As the last case study, market equilibrium price shows that none of the numerical ODE solvers can estimate its behavior which is due to its sudden changing nature.

## 1. Introduction

The utilization of Ordinary Differential Equation (ODE) solvers is fundamental in various scientific and engineering disciplines, providing powerful tools for modeling and analyzing complex dynamic systems [1]. More recently, ODE-based methods have gained prominence in electro/mechanical engineering [2-7], artificial intelligence (8-10), machine learning (11,12), optimization (13-17), cybersecurity (18-20) and neural network research (21,22). Applications include continuous-depth models like Neural ODEs, dynamic system modeling in reinforcement learning, and the formulation of optimization algorithms as continuous-time processes, offering new perspectives on learning dynamics and model interpretability (23). In mechanical engineering, ODE solvers have become indispensable for modeling intricate mechanical systems, including vibrations and multi-body dynamics, with methods like Euler's method, Hyun's method, and the midpoint method serving as efficient and accurate simulation tools. In the realm of electrical engineering, ODE solvers find applications in circuit analysis and control system design, where techniques like the fourth-order Runge-Kutta (RK4) method are instrumental in capturing transient responses and ensuring system stability (24).

The choice of an appropriate ODE solver is crucial, as numerous programming environments offer a range of solvers, each equipped with distinct numerical techniques for integrating systems of ODEs. The decision between explicit and implicit methods significantly influences the efficacy of ODE solutions, guided by the specific solver selected (25).

ODE solvers play a pivotal role in modeling population dynamics, with previous studies relying on various versions of the exponential model to project changes in population change (26). The adoption of the logistic model enhances precision in these projections, as it incorporates multiple data types (27). In a different context, ODE solvers are invaluable for predicting changes in the internal temperature of buildings. While prior research has extensively addressed the issue of thermal comfort in educational and residential buildings, a primary focus on specific educational facilities has been evident (28-31). Moreover, many of these studies have failed to distinguish between thermal comfort during occupancy and non-occupancy periods, primarily monitoring buildings during working hours (32). In the third case study, our attention shifts to determining the equilibrium price in the market, aiming to establish more precise relationships that accurately reflect the reference data under consideration.

Differential equations form the mathematical foundation for describing dynamic change, evolution, and variation in various contexts within mathematics and the natural and social sciences. These equations connect differentials, derivatives, and functions to describe dynamic phenomena, often defining quantities as the rate of change of other quantities or gradients, which are integral to the formulation of differential equations (33).

In the domain of numerical ODE solvers, a range of methods plays pivotal roles in engineering applications. Euler's method, known for its simplicity and rapid calculations, coexists with Heun's method, which enhances accuracy through dynamic step size adaptation (34). The midpoint method strikes a balance between accuracy and computational efficiency. Among these techniques, the RK4 method stands out, offering higher-order accuracy and gaining prominence across multiple engineering disciplines. Collectively, these solver techniques empower engineers to model intricate real-world systems, ensuring reliable predictions and informed decision-making for a broad spectrum of engineering challenges (35).

The Ordinary Differential Equation (ODE) approach finds extensive utilization across diverse scientific and engineering domains, encompassing fields such as biology, economics, and physics. It serves as a valuable tool for comprehensively analyzing the behavior of intricate systems (23).

In mechanical engineering, researchers have employed ODE solvers to model complex mechanical systems, including vibrations and multi-body dynamics, utilizing methods like Euler's method, Hyun's method, and the midpoint method for efficient and accurate simulations. In electrical engineering, ODE solvers have found application in circuit analysis and control system design, harnessing techniques like the fourth-order Runge-Kutta (RK4) method to capture transient responses and system stability (24).

In this study we examine the accuracy and robustness of 5 different ODE solvers namely: Euler's, Heun's, Midpoint, Runge-kutta 4th order and ODE45 on three different case studies namely: population dynamics, temperature change and market equilibrium price. The outcome of this study will shed light on each of the mentioned ODE solvers leading to an improved understanding on their application in future physical case studies.

## 2. Methodology

Without computers, ODEs are usually solved with analytical integration techniques. The important fact is that exact solutions for many ODEs used in physical applications are not available. Also, numerical methods offer the only viable alternative for these cases. Because these numerical methods usually require computers, engineers in the pre-computer era were somewhat limited in the scope of their investigations. All the ode solvers follow a general formula of the following form

$$\frac{d_y}{d_x} = f(x, y) \tag{1}$$

In which all the new values for the next time step can be generated using: New value = old value + slope × step size or, in mathematical terms,

$$y_{i+1} = y_i + \Phi h$$

All one-step methods can be expressed in this general form, with the only difference being the manner in which the slope is estimated. With that said we use five different ODE solvers Namely Euler's, Heun's, Midpoint and 4$^{th}$ order Runge Kutta Methods along with ODE45.

### 2.1 Euler's Method

In this method, a direct estimate of the slope at $x_i$ will be provided by the first derivative as

$$\Phi = f(x_i, y_i) \tag{2}$$

where $f(x_i, y_i)$ represents the differential equation evaluated at $x_i$ and $y_i$. This estimate can be used together with

$$y_{i+1} = y_i + f(x_i, y_i)h \tag{3}$$

This formula is referred to as Euler's (also called as the Euler-Cauchy or the point-slope) method. A new value of y is generated using the slope (the value of the first derivative at the original value of x) to linearly extrapolate over the step size h.

The numerical solution of ODEs will be evaluated by two different error types.

1. Truncation, or discretization, errors generated by the nature of the techniques used to estimate values of y.
2. Round-off errors generated by the limited numbers of significant digits that can be calculated by a computer.

The truncation errors are composed of two terms. The first term is a local error that is caused by the application of the method in use over a single step. The second term is a propagated truncation error that is caused by the approximations generated during the previous steps. The sum of the two terms composes the total, or global, truncation error. Insight into the value and properties of the truncation error can be achieved by deriving Euler's method directly from the Taylor series expansion. Notably, the differential equation being integrated will be of the general form

$$y' = f(x, y) \tag{4}$$

where $y' = {d_y}/{d_x}$ and x and y are the independent and the dependent variables, respectively. If the solution—that is, the function describing the behavior of y—has continuous derivatives, it can be represented by a Taylor series expansion about a value ($x_i$, $y_i$) as in

$$y_{i+1} = y_i + y'h + \frac{y_i^{''}}{2!} + \cdots + \frac{{y_i}^{(n)}}{n!}h^n + R_n \tag{5}$$

where h is $x_i$+1 − $x_i$ and $R_n$ is the remainder term, defined as

$$R_n = \frac{y^{(n+1)}(\xi)}{(n+1)!}h^{n+1} \tag{6}$$

where ξ lies somewhere in the interval from $x_i$ to $x_{i+1}$. An alternative form can be developed by substituting equation (6) into equation (5) and (4) to yield

$$y_{i+1} = y_i + f(x_i, y_i)h + \frac{f'(x_i, y_i)}{2!}h^2 + \cdots + \frac{f^{(n-1)}(x_i, y_i)}{n!}h^n + O(h^{n+1}) \tag{7}$$

where $O(h^{n+1})$ specifies that the local truncation error is proportional to the step size raised to the $(n + 1)^{th}$ power. By comparing equation (3) and (7), Euler's method corresponds to the Taylor series up to and including the term $f(x_i, y_i)h$. Additionally, the comparison indicates that a truncation error occurs because we approximate the true solution using a finite number of terms from the Taylor series. We thus truncate, or leave out, a part of the true solution. For example, the truncation error in Euler's method is attributable to the remaining terms in the Taylor series expansion that were not included in equation (3). Subtracting equation (3) from equation (7) yields

$$E_t = \frac{f'(x_i, y_i)}{2!}h^2 + \cdots + O(h^{n+1}) \tag{8}$$

where Et is equal to the true local truncation error. For sufficiently small h, the errors in the terms in equation (8) usually decrease as the order increases, and the result is often represented as

$E_a = \frac{f'(x_i,y_i)}{2!}h^2$ Or $E_a = O(h^2)$

where $E_a$ is the approximate local truncation error.

## 2.2 Heun's Method

An effective source of error in Euler's method is that the derivative at the beginning of the interval is assumed to be valid across the entire interval. Two types of improvements are available to help overcome this shortcoming. Both modifications belong to a larger class of solution techniques called Runge-Kutta methods. However, due to having a very straightforward graphical interpretation, we need to present them prior to their formal derivation as Runge-Kutta methods.

One of the methods which improves the evaluation of the slope involves the determination of two derivatives for the interval—the first one at the initial point and the second at the end point. The two derivatives are then averaged to obtain an improved estimate of the slope for the entire interval. Please note that in Euler's method, the slope at the beginning of an interval as below is used to extrapolate linearly to $y_{i+1}$:

$$y'_i = f(x_i, y_i) \tag{9}$$

For the Euler method, we ought to stop at this point. However, in Heun's method the $y^0_{i+1}$ calculated in equation (9) is not the result, but an intermediate prediction. This is why we have distinguished it with a superscript 0. Equation (9) is called a predictor equation. It provides an estimate of $y_{i+1}$ that allows the evaluation of a calculated slope at the end of the interval:

$$y'_{i+1} = f(x_{i+1}, y^0_{i+1}) \tag{10}$$

Thus, the two slopes, in equation (9) and (10), can be combined to obtain an average slope for the interval:

$$\overline{y'} = \frac{y'_i + y'_{i+1}}{2} = \frac{f(x_{i+1}, y^0_{i+1})}{2} \tag{11}$$

This average slope is then used to extrapolate linearly from $y_i$ to $y_{i+1}$ using Euler's method:

$$y_{i+1} = y_i + \frac{f(x_i,y_i) + f(x_{i+1},y^0_{i+1})}{2}h. \tag{12}$$

This is called a corrector equation.

The Heun method is a predictor-corrector approach. All the multistep methods are of this type. As derived above, it can be expressed as

Predictor $$y^0_{i+1} = y_i + f(x_i, y_i)h \tag{13}$$

Corrector $$y_{i+1} = y_i + \frac{f(x_i,y_i) + f(x_{i+1},y^0_{i+1})}{2}h \tag{14}$$

Note that because equation (14) has $y_{i+1}$ on both sides of the equation, it can be considered in an iterative fashion. That is, an old value can be used repeatedly to provide an improved estimate of $y_{i+1}$. It should be noted that this iterative process does not converge on the true answer yet will converge on an estimate with a finite truncation error. As with similar iterative methods, a termination criterion for convergence of the corrector is shown by

$$|\varepsilon_a| = \left|\frac{y_{i+1}^{j} - y_{i+1}^{j-1}}{y_{i+1}^{j}}\right| 100\% \qquad (15)$$

where $y_{i+1}^{j-1}$ and $y_{i+1}^{j}$ are the result from the prior and the present iteration of the corrector, respectively.

## 2.3 Midpoint Method

This technique utilizes Euler's method to estimate a value of $y$ at the midpoint of the interval

$$y_{i+1/2} = y_i + f(t_i, y_i)\frac{h}{2}\% \qquad (16)$$

Then, this predicted value is used to calculate a slope at the midpoint:

$$y'_{i+1/2} = f(t_{i+1/2}, y_{i+1/2})\ \% \qquad (17)$$

which is considered to represent a valid approximation of the average slope for the entire interval. This slope is then used to extrapolate linearly from $t_i$ to $t_{i+1}$

$$y_{i+1} = y_i + f\left(t_{i+1/2}, y_{i+1/2}\right)h\ \% \qquad (18)$$

Observe that because $y_{i+1}$ is not on both sides, the corrector cannot be applied iteratively to improve the solution as was done with Heun's method. Similar to Heun's method, the midpoint method can also be linked to Newton-Cotes integration formulas (36).

## 2.4 Runge-Kutta Method

Runge-Kutta (RK) methods achieve the accuracy of a Taylor series approach without the need for the evaluation of higher derivatives. Many variations exist but all can be shown in the generalized form as below

$$y_{i+1} = y_i + \phi h \qquad (19)$$

where $\varphi$ is called an *increment function,* which can be interpreted as a representative slope $i$+1⁄2 over the interval. The increment function can be written in general form as

$$\phi = a_1 k_1 + a_2 k_2 + \cdots + a_n k_n \tag{20}$$

where the $a$'s are constants and the $k$'s are

$$k_1 = f(t_i, y_i) \tag{21}$$
$$k_2 = f(t_i + p_1 h, y_i + q_{11} k_1 h)$$
$$k_3 = f(t_i + p_2 h, y_i + q_{21} k_1 h + q_{22} k_2 h)$$
$$\vdots$$
$$k_n = f(t_i + p_{n-1} h y_i + q_{n-1,1} k_1 h + q_{n-1,2} k_2 h + \cdots + q_{n-1,n-1} k_{n-1} h$$

where the $p$'s and $q$'s are constants. Bear in mind that the $k$'s are recurrence relationships. That is, $k_1$ appears in the equation for $k_2$ which appears in the equation for $k_3$ and like the rest. Because each $k$ is a functional evaluation, this recurrence makes RK methods efficient for computer programs. By setting appropriate values of a we can derive all the ODE solvers described previously from the generalized formula.

Various types of RK methods can be devised by employing different types of terms in the increment function as specified by $n$. Note that the first-order RK method with $n$=1 is Euler's method. At least for the lower-order versions, the number of terms $n$ usually represents the order of the approach.

## 2.5 Runge-Kutta 4th Order Method

The most versatile RK methods are 4th order. As with the second-order approaches, there are an infinite number of types. The following is the most popular used form, and we therefore call it the *classical fourth-order RK method:*

$$y_{i+1} = y_i + \frac{1}{6}(k_1 + 2k_2 + 2k_3 + k_4)h \tag{22}$$

where:

$$k_1 = f(t_i, y_i)$$

$$k_2 = f(t_i + \frac{1}{2}h, y_i + \frac{1}{2}k_1 h)$$

$$k_3 = f(t_i + \frac{1}{2}h, y_i + \frac{1}{2}k_2 h)$$

$k_4 = f(t_i + h, y_i + k_3 h)$

Notice that for ODEs that are a function of the independent variable alone, the classical 4$^{th}$-order RK method is similar to Simpson's 1⁄3 rule. In addition, the 4$^{th}$-order RK method is similar to the Heun's approach in that multiple evaluation of the slope are developed to come up with an improved average slope for the interval.

### 2.6 ODE45 Method

ODE45 is a numerical solver for ordinary differential equations (ODEs) available in MATLAB, and its methodology involves implementing a variable-step, variable-order algorithm. The ode45 function in MATLAB is a versatile numerical ODE solver that utilizes a variable-step, variable-order integration method to numerically solve ordinary differential equations (ODEs) of the form $y' = {d_y}/{d_x}$ where y is the vector of dependent variables and x is the independent variable. By providing the ODE function f, initial conditions, and the desired range of x values, ode45 dynamically adjusts its integration step size based on the ODE behavior, using a combination of lower-order and higher-order methods to ensure accuracy and efficiency. When provided with an initial value problem, ODE45 employs a combination of explicit Runge-Kutta (4th and 5th order) formulas to approximate the solution. It starts by taking an initial step and monitors the error between the 4th and 5th order approximations. If the error is within the specified tolerance, the step size is adjusted, and the process continues. If the error exceeds the tolerance, the solver reduces the step size to ensure accuracy. ODE45 dynamically adapts the step size throughout the integration, which allows it to efficiently handle both stiff and non-stiff ODEs. This method provides a robust and versatile approach to solving a wide range of ODEs encountered in various scientific and engineering applications while maintaining a balance between accuracy and computational efficiency.

For each case study error has been calculated for each numerical differentiation method versus both empirical relationship and experimental. Consider that the data resulted from empirical equation can be shown as $E_1$, $E_2$, $E_3$… $E_n$, and the data from the experimental work is shown by $P_1$, $P_2$, $P_3$… $P_n$ as well as the data form each numerical differentiation method as $X_1$, $X_2$, $X_3$… $X_n$ we can calculate the error for each method as

$$Error_{wrt_Exp} = \frac{\sum_{i=1}^{n}(P_i - X_i)}{\sum_{i=1}^{n}(P_i)} \quad (23)$$

$$Error_{wrt_Emp} = \frac{\sum_{i=1}^{n}(E_i - X_i)}{\sum_{i=1}^{n}(E_i)} \quad (24)$$

## 3. Results and Discussion

The section delves into the comprehensive evaluation and comparative analysis of various ordinary differential equation solvers across three distinct case studies. It provides a detailed graphical view of the numerical results obtained from the application of these solvers to the selected case studies, offering insights into their respective performances and computational efficiency. The discussion unfolds through an in-depth examination of the implications of these findings, exploring how solver selection impacts the accuracy and stability of solutions in real-world scenarios. Moreover, it scrutinizes the computational resources required by each solver and their suitability for different types of ODE problems, shedding light on the practical implications for researchers working in diverse scientific and engineering domains.

### 3.1 Population Dynamics

In the figure below, you clearly see the comparison of different ODE solvers on population change. It is illustrated that for this case the results of different solvers are very close with each other so that we can say that the best method is the one with lowest solution cost. Naturally, Euler’s method has the lowest solution cost in case of numerical analysis and after that Heun’s and Midpoint Method are of the same order. As the greatest solution cost among all these numerical ODE solvers, RK $4^{th}$ order can be considered as the last method. The error presented below indicates the difference between the types of the solver and empirical/experimental data. Regardless of the solution cost of the solvers form figure 1, we can conclude that the lowest error is associated with error of Heun’s solver, while the highest error is linked to Euler’s solver with respect to experimental data.

When assessing the exponential model from table 1, it is evident that Euler's method exhibited the highest error, with experimental error measuring 2.9281, signifying a relatively larger deviation from the expected values. In contrast, the Heun’s method demonstrated slightly better accuracy, as indicated by experimental error at 2.9438. The midpoint method and the RK method performed with similar accuracy, yielding experimental errors at 2.9418 and 2.9417, respectively. However, the ode45 solver stood out with the lowest error at 2.9417, signifying its superior accuracy in approximating the population dynamics using the exponential model. This analysis highlights the varying levels of precision among these solvers, with ode45 delivering the most reliable results in this specific context.

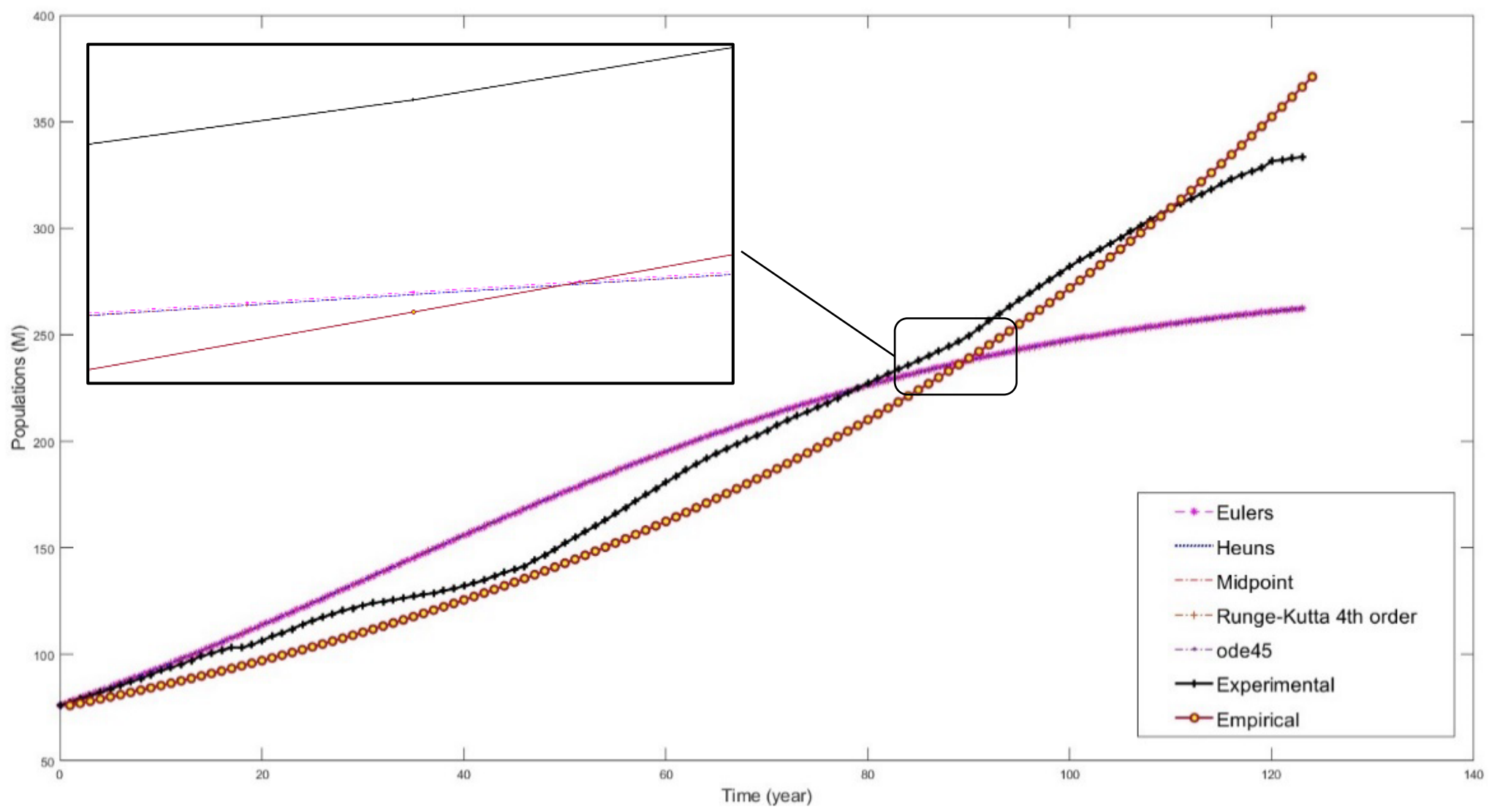

Figure 1. Comparison of 5 considered ODE solvers on population dynamics

Table 1. Error values of considered ODE solvers on population dynamics with respect to both experimental and empirical results

| Error | value |
|---|---|
| $e_{exp}^{euler}$ | 2.9281 |
| $e_{emp}^{euler}$ | 0.2020 |
| $e_{exp}^{heun}$ | 2.9438 |
| $e_{emp}^{heun}$ | 0.1858 |
| $e_{exp}^{midpoint}$ | 2.9418 |
| $e_{emp}^{midpoint}$ | 0.1878 |
| $e_{exp}^{R-K}$ | 2.9417 |
| $e_{emp}^{R-K}$ | 0.1879 |
| $e_{exp}^{ode45}$ | 2.9417 |
| $e_{emp}^{ode45}$ | 0.1879 |

### 3.2 Temperature Change

As a different case study having a new nature of changing, we see in picture below that the behavior of the ODE solvers in estimating the results are different from the previous method. We see that the results of all solvers are close to each other except in the initial interval of the curve where Euler's method shows a behavior relatively off the empirical results. Regardless of the solution cost of the solvers form 2, we can conclude that the lowest error is associated with Midpoint solver with respect to experimental data, while the highest error is linked to Euler's solver with respect to empirical data. Therefore, the best method for estimating and later predicting the results of temperature change in a building is Midpoint method which is an interesting result while we have the presence of other methods with higher solution cost and terms of calculation. This phenomenon shows the importance of the nature of each model for choosing the best ODE solver.

When considering the results of table 2, Euler's method demonstrated an experimental error at 0.8122, indicating a moderate level of precision in approximating temperature changes. The Heun's method exhibited improved accuracy with an error at 0.5367, suggesting a closer alignment with the expected temperature dynamics. Similarly, the midpoint method presented a similar level of accuracy, with an error at 0.5362, offering reliable temperature predictions. The RK method displayed a slightly higher error at 0.6182, indicating a marginally lower precision in modeling temperature changes. Remarkably, the ode45 solver proved to be the most accurate among the solvers, with the lowest error at 0.6790, signifying its superior performance in approximating temperature variations using the exponential model. This analysis underscores the varying levels of accuracy among the solvers, with ode45 emerging as the most reliable choice for this particular temperature change case study.

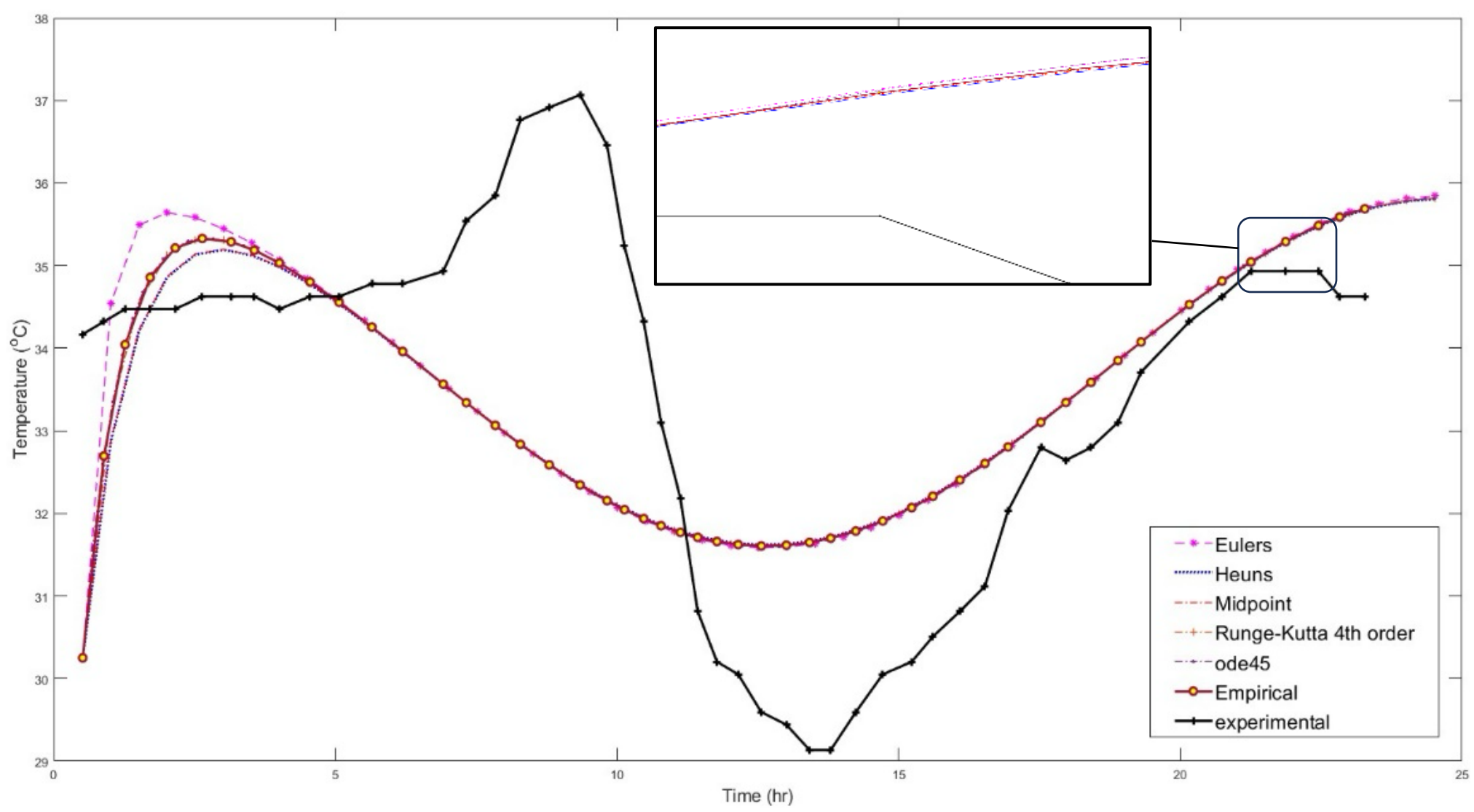

Figure 2. Comparison of 5 considered ODE solvers on temperature change

Table 2. Error values of considered ODE solvers on temperature change with respect to both experimental and empirical results

| Error | value |
|---|---|
| $e_{exp}^{euler}$ | 0.8122 |
| $e_{emp}^{euler}$ | 1.0929 |
| $e_{exp}^{heun}$ | 0.5367 |
| $e_{emp}^{heun}$ | 0.8166 |
| $e_{exp}^{midpoint}$ | 0.5362 |
| $e_{emp}^{midpoint}$ | 0.8161 |
| $e_{exp}^{R-K}$ | 0.6182 |
| $e_{emp}^{R-K}$ | 0.8983 |
| $e_{exp}^{ode45}$ | 0.6790 |
| $e_{emp}^{ode45}$ | 0.8891 |

### 3.3 Market Equilibrium Price

In figure below, we can compare the result of estimation of 5 ODE solver with both empirical and experimental data for market equilibrium price. We see from figure 3 that for this model all the ODE solvers after a short period of time fail to estimate the results and blow up in prediction. This is because of the nature of this model and due to the steep slope of the model. This behavior makes the denominator of all the solvers to be a very small number and thus causes this failure of prediction.

Regardless of the inability of estimating the price change for all the methods from table 3 we see that Euler's method exhibited a relatively high error at 1.2896e+03, suggesting a notable deviation from the expected market equilibrium price. The Heun's method showed a similar level of inaccuracy, with an error at 1.0228e+04, indicating a substantial margin of error. Likewise, the midpoint method exhibited a comparable level of imprecision, with an error at 1.0228e+04, further highlighting the inability of these methods to accurately model market equilibrium price dynamics. The RK method displayed the highest error among the solvers at 2.1990e+04, signifying a considerable discrepancy from the expected values. Notably, the ode45 solver also presented a high error at 2.2770e+04, indicating a significant margin of inaccuracy in approximating market equilibrium prices using the exponential model. This analysis underscores the limitations of these solvers in accurately modeling market equilibrium prices, suggesting the need for alternative approaches or models for improved predictions.

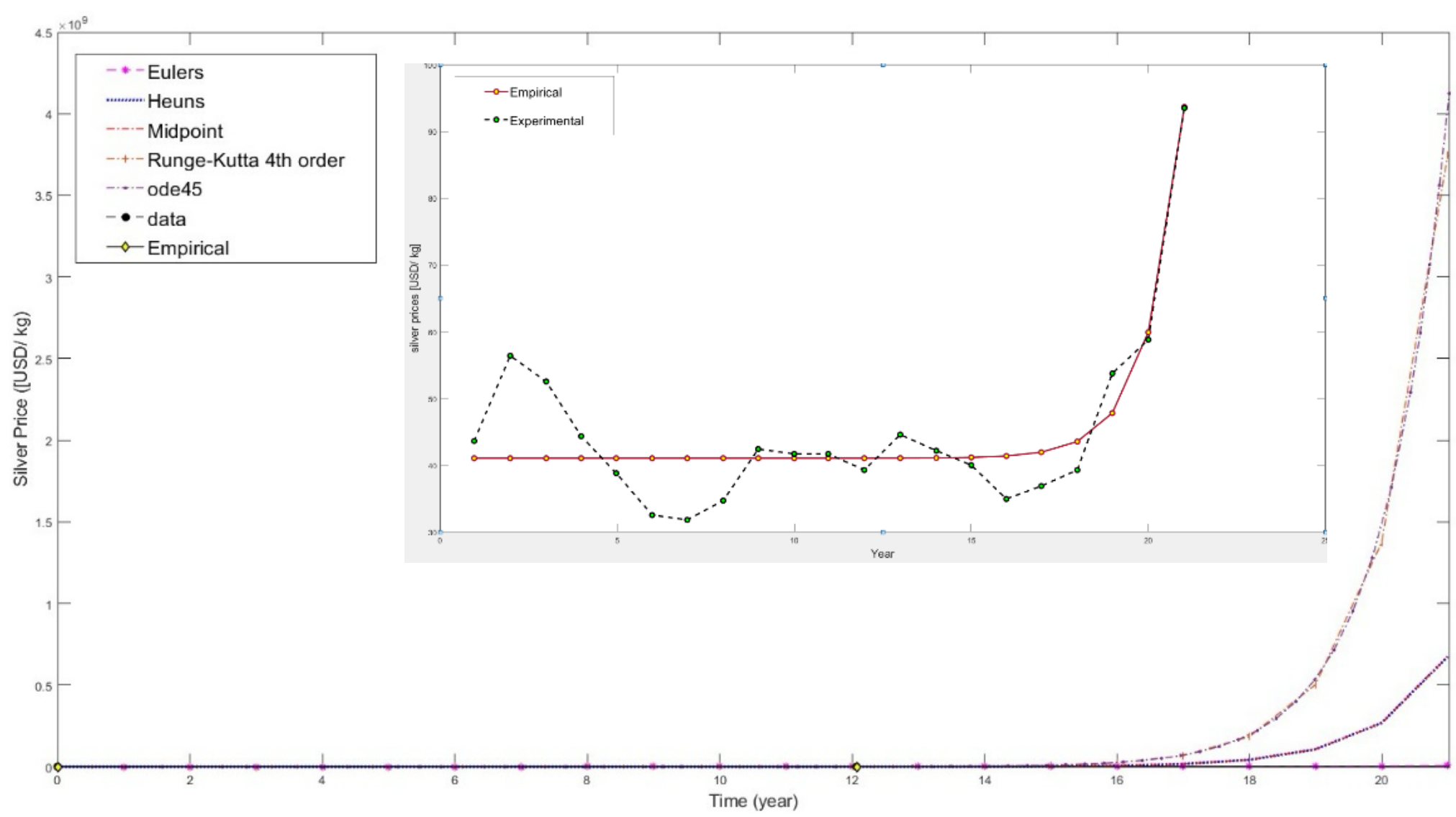

Figure 3. Comparison of 5 considered ODE solvers on market equilibrium price

Table 3, Error values of considered ODE solvers on market equilibrium price with respect to both experimental and empirical results

| Error | Value |
|---|---|
| $e_{exp}^{euler}$ | 1.2896e+03 |
| $e_{emp}^{euler}$ | 100 |
| $e_{exp}^{heun}$ | 1.0228e+04 |
| $e_{emp}^{heun}$ | 100 |
| $e_{exp}^{midpoint}$ | 1.0228e+04 |
| $e_{emp}^{midpoint}$ | 100 |
| $e_{exp}^{R-K}$ | 2.1990e+04 |
| $e_{emp}^{R-K}$ | 100 |
| $e_{exp}^{ode45}$ | 2.2770e+04 |
| $e_{emp}^{ode45}$ | 100 |

## 4. Conclusion

In this work we consider three different case studies namely: logistic model for population change, temperature change of a building and market equilibrium price to reveal the accuracy of estimation of results of 5 ODE solvers. The considered ODE solvers are Euler's, Heun's, Midpoint, Runge-kutta 4$^{th}$ order and ODE45 method. We also show the results of previous experimental works to evaluate the results shown in this work which shows that they are in a good agreement. The results illustrate that for the first case study all of the numerical solvers show similar accuracy, and we can just consider the method with smaller solution cost as the best solver. For the second case study, we see that in some intervals of time due to the behavior of the empirical model the solvers deviate from each other and considering the relative error of the solvers it turns out that Midpoint method is the best solver. Finally, for the last case study, we see that none of the numerical solvers can estimate (and for future works predict) the behavior of the model. This phenomenon is due to the sudden change in the behavior of the model. As a general result we can assert that based the behavior of each model the application and accuracy of each solver is different.